\documentclass{amsart}
\usepackage{amsthm,amsmath,amsfonts,amssymb,amscd,mathrsfs,graphics}
\usepackage{latexsym}
\usepackage{graphicx}

\usepackage{epic}

\usepackage{enumerate}
\usepackage[...]{youngtab}

\usepackage{tikz}

\newtheorem{thm}{Theorem}[section]

\newtheorem{lem}[thm]{Lemma}
\newtheorem{prop}[thm]{Proposition}
\newtheorem{cor}[thm]{Corollary}

\makeatletter
\renewcommand{\@seccntformat}[1]{\S{\csname
the#1\endcsname}\hspace{0.5em}}
\makeatother

\begin{document}

\title{Difference sets disjoint from a subgroup II: groups of order $4p^2$}

\author{ Stephen P.  Humphries, Nathan L. Nicholson}
  \address{Department of Mathematics,  Brigham Young University, Provo, 
UT 84602, U.S.A.
E-mail: steve@mathematics.byu.edu, nlnicholson24@gmail.com}
\date{}
\maketitle

\begin{abstract}  
We study finite groups $G$ having a normal subgroup $H$ and $D \subset G \setminus H, D \cap D^{-1}=\emptyset,$ such that the multiset $\{ xy^{-1}:x,y \in D\}$ has every non-identity element occur the same number of times (such a $D$ is called a {\it DRAD difference set}). We show that there are no such groups of order $4p^2$, where $p$ is an odd  prime.\medskip


\medskip
\end{abstract}

\theoremstyle{plain}

\theoremstyle{definition}
\newtheorem*{dfn}{Definition}
\newtheorem{exa}[thm]{Example}
\newtheorem{rem}[thm]{Remark}

\newcommand{\ds}{\displaystyle}
\newcommand{\bs}{\boldsymbol}
\newcommand{\mb}{\mathbb}
\newcommand{\mc}{\mathcal}
\newcommand{\mf}{\mathfrak}
\renewcommand{\mod}{\operatorname{mod}}
\newcommand{\mult}{\operatorname{Mult}}

\def \a{\alpha} \def \b{\beta} \def \d{\delta} \def \e{\varepsilon} \def \g{\gamma} \def \k{\kappa} \def \l{\lambda} \def \s{\sigma} \def \t{\theta} \def \z{\zeta}

\numberwithin{equation}{section}

\setlength{\leftmargini}{1.em} \setlength{\leftmarginii}{1.em}
\renewcommand{\labelenumi}{\setlength{\labelwidth}{\leftmargin}
   \addtolength{\labelwidth}{-\labelsep}
   \hbox to \labelwidth{\theenumi.\hfill}}

\maketitle

\section{Introduction}

 For a group $G$ we will identify a finite subset $X \subseteq G$ with the element $\sum_{x \in X} x \in \mathbb Q G$ of the group algebra. We also let $X^{-1}=\{x^{-1}:x \in X\}$.  Write $\mathcal C_n$ for the cyclic group of order $n$.

 A $(v,k,\lambda)$ {\it difference set} is a subset $D \subset G, |D|=k$,  such that every element $1 \ne g \in G$ occurs $\lambda $ times in the multiset $\{xy^{-1}:x,y \in D\}$. Here  $|G|=v$.
 
Then a $(v,k,\lambda)$ difference set is a $(v,k,\lambda)$ {\it DRAD difference set (with subgroup $H$ and difference set $D$)} if it also satisfies the   conditions:
 there is a subgroup $1\ne H \triangleleft  G$ such that

\noindent (1) $D \cap D^{-1}=\emptyset$;

\noindent (2) $G \setminus (D \cup D^{-1})=H $.

A group $G$ will be called a {\it DRAD difference set group} if there is a DRAD difference set over $G$. 
  See   \cite {dp,ito1,ito2} for more on DRADs. 
  DRAD difference sets are examples of Hadamard (or Menon) difference sets \cite {DJ}. 
  Let $$ h=|H|,\quad u=|G:H|.$$

We proved the following result in \cite {DS1}:
\begin{thm} \label{thm1}  Let $G$ be a $(v,k,\lambda)$  DRAD difference set group with subgroup $H$ and difference set $D$. Then 

\noindent (i) $  u=h\ge 4$ is even, $v=|G|=h^2$, and  
$$\lambda=\frac 1 4  {h(h-2)},\,\,\,  k=\frac 1 2 h(h-1);$$

\noindent (ii)  each non-trivial coset  $Hg\ne H$ meets $D$ in $h/2$ points;

\noindent (iii) $H$ contains the subgroup generated by all the involutions in $G$;

\noindent (iv) any abelian $(v,k,\lambda)$ DRAD difference set group is a $2$-group.
\end{thm}

All known DRAD difference set groups are $2$-groups. There is one such group of order $16$ and at least $33$ of order $64$. In \cite {DS1} a bi-infinite family of DRAD groups is constructed.
In this paper we show 
\begin{thm} \label{thm2}  
There are no $(v,k,\lambda)$ DRAD difference set groups of order  $4p^2$, for an odd prime $p$.
\end{thm}

Theorem \ref {thm2}  is related to Question 9 of \cite {DJ}, where Davis and Jedwab ask if there are any Hadamard difference sets of order $4p^2$, where $p>5$ is an odd prime. Note that when $h=2p$, then by Theorem \ref {thm1} we have $(v,k,\lambda)=(4p^2,2p^2-p,p^2-p)$, the parameters for a Hadamard difference set.

For the the proof of Theorem \ref {thm2}  we make use of a result of Iiams \cite {ii}, who
 showed that any group of order $4p^2$ (where $ p>3$ is a prime) that has $\mathcal C_p \times \mathcal C_2^2$ as a factor group, does not have a $(4p^2,2p^2-p,p^2-p)$ difference set. See Theorem \ref {thmii} below for a  full statement of the  result of Iiams.
This result was extended by Wan \cite {wan}   to show that any group of order $4p^4$ that has 
$\mathcal C_{p^2} \times \mathcal C_2^2$ as a factor group 
does not have a $(4p^4,2p^4-p^2,p^4-p^2)$ difference set.
Lastly, the authors of \cite {as}  show that if $q=p^n, p>3$ a prime, where $G$ has $\mathcal C_q \times \mathcal C_2^2$ as a factor group, then $G$ does not admit a $(4q^2,2q^2-q,q^2-q)$ difference set.
We note that Smith \cite {sm} was the first to find non-abelian examples with these parameters, and that McFarland \cite{mcf} shows 
that no nontrivial difference set exists in abelian groups of order $4p^2$ where $ p>3$ is a prime. We also note that Davis and Polhill construct an infinite family for abelian groups \cite {dp}.

\noindent {\bf Acknowledgement} We thank a referee for pointing out an error in an earlier version of this paper. All calculations made in writing this paper were accomplished using Magma \cite {Ma}.  

\section{The  result of Iiams and the cases $G_4, G_{13}, G_{16}$}

The result of Iiams \cite {ii} that we use is
\begin{thm} \label{thmii} Let $p\ge 5$ be prime, and let $G$ be a group of order $4p^2$ that contains a Menon-type difference set. Then either

\noindent (i) $p \equiv 1 \mod 4$ and $G\cong G_{11}$;

\noindent (ii) $G$ is isomorphic to one of $G_4,G_{13},G_{14},G_{15}, G_{16}$.
\end{thm}

Here the groups $G_{11},G_4,G_{13},G_{14},G_{15}, G_{16}$ are defined below.
We now consider these six cases. (In general we have replaced $z$ by $z^{-1}$ in the presentations that Iiams gives, thus allowing   conjugation to be written $y^z=z^{-1}yz$.)

\medskip 

\noindent {\bf Case $G_4$}: Let $p$ be an odd  prime and let  $f \in \mathbb N$ such that $f^2 \equiv -1 \mod p^2$. Then 
$G_4$ has the following presentation:
$$G_4=\langle x,z|x^{p^2}, z^4,x^z=x^f\rangle.$$

Note that $x^{z^2}=x^{-1}$, which gives $(z^2)^x=x^{-2}z^2,$ and so
$(z^2)^{x^k}=x^{-2k}z^2.$

Thus  each  $x^{-2k}z^2$  has order $2$, and since $x$ has order ${p^2}$ (with $p$ odd) we see that all these elements of order $2$ generate the subgroup
$\langle x,z^2\rangle$ of order $|G_4|/2=h^2/2$. By Theorem \ref{thm1} we have   $|H|=h$ and $H$ has to contain this subgroup,  so we have $h^2/2 \le |H|=h$,  a contradiction to Theorem \ref {thm1} (i). 
\medskip

\noindent {\bf Case $G_{13}$}:  Here $G_{13}$ has presentation:
$$G_{13}=\langle x,y,z|x^p,y^p,(x,y),z^4,x^z=y^{-1},y^z=x\rangle.$$
Again we have $x^{z^2}=x^{-1},y^{z^2}=y^{-1}$, which gives
$(z^2x^ay^b)^2=1$ for $0\le a,b<p$.
 Thus, as in the $G_{4}$ case,  the subgroup $\langle z^2x^ay^b:0\le a,b <p\rangle \le G_{13}$ generated by the involutions,
has index $2$. This rules out this case.

\medskip

\noindent {\bf Case $G_{16}$}:  Let   $f \in \mathbb N$ such that $f^2 \equiv -1 \mod p$. Then  
$$G_{16}=\langle x,y,z|x^p,y^p,(x,y),z^4,x^z=x^f,y^z=y^f\rangle.$$
Again we have $x^{z^2}=x^{-1},y^{z^2}=y^{-1}$, and so, as in the $G_{13}$ case,  the subgroup of $G_{16}$ generated by the involutions
has index $2$, and we obtain a contradiction. 
\medskip

\section {The groups $G_{11}$ and $G_{14}$}

For these cases we will need the following Lemma:
\begin{lem} \label{lemwww} Suppose that $G$ has a non-principal linear character $\chi$. 

  If $\chi(H)=0$ and $\chi$  takes values in a field $K$ where $i=\sqrt{-1}\notin K,$ then $G$ is not a DRAD group with subgroup $H$.

\end{lem}
\noindent {\it Proof} 
Since $\chi$ is linear and non-principal we have $\chi(G)=0.$
By hypothesis we have $\chi(H)=0$, so that $G=H+D+D^{-1}$ gives $\chi(D^{-1})=-\chi(D).$ Then 
$$
DD^{-1}=\lambda G+(k-\lambda)=\lambda G+\frac 1 4 h^2
$$
 and the linearity of $\chi$ gives $\chi(D)^2=-\frac 1 4 h^2$, which gives $\chi(D)=\pm i \frac h 2\in K$, a contradiction.\qed  
\medskip

\noindent {\bf Case $G_{11}$}:  Let    $p$ be a prime where $p \equiv 1 \mod 4$. Then 
$$G=G_{11}=\langle x,y,z|x^p,y^p,(x,y),z^4,x^z=x,y^z=y^{-1}\rangle\cong \langle x \rangle  \times \langle y,z\rangle.$$ 
One checks that there are two possibilities for $H$:
$Z=Z(G)=\langle x,z^2\rangle$ and $Y=\langle y,z^2\rangle.$ Note that   $Z,Y \triangleleft G$.

Let $\zeta_p=\exp{2\pi i/p}.$ 

First consider  $H=Z=\langle x,z^2\rangle$ and let $\chi$ be the linear character 
\begin{align}
\label{eqn0}
\chi(x)=\zeta_{p},\quad \chi(y)=1, \quad \chi(z)=-1.
\end{align}
Then $\chi(H)=0$ (since $x \in H$ and $\chi(x) \ne 1$) and $\chi(g) \in \mathbb Q(\zeta_{p})$ for all $g \in G$,
from which   Lemma \ref {lemwww}   shows that this case cannot happen, since $i \notin \mathbb Q(\zeta_p).$
\medskip 

So now consider  $H=Y=\langle y,z^2\rangle\triangleleft G.$

Elements of $G$ will have {\it normal form} $g_{j,k,w}=x^jy^kz^w, 0\le j,k<p,0\le w<4$. Here elements of $H$ have the form $g_{0,k,w}=y^kz^w, 0\le k<p,w \in \{0,2\}$.
Now we  let 
$$D=\sum_{j,k,w} \varepsilon_{j,k,w}g_{j,k,w},$$ where  $\varepsilon_{j,k,w} \in \{0,1\}$, and   $\sum_{j,k,w}$ means   we sum over 
         $0\le j,k<p,w \in \{0,1,2,3\}.$

For some fixed $m \in \mathbb Z$ we consider the linear character 
\begin{align}
\label{eqn000}
 \chi_m: x\mapsto \zeta_p^m,\quad y \mapsto 1,\quad z \mapsto i=\sqrt{-1}.
\end{align}
Then we have $\chi_m(H)=0$, since $z^2 \in H$ and $\chi_m(z^2)=-1$. We also have $\chi_m(G)=0$. Thus  $D+D^{-1}=G-H$ gives $\chi_m(D^{-1})=-\chi_m(D)$
and so $DD^{-1}=\lambda G+(k-\lambda)$ gives 
$$-\chi_m(D)^2=k-\lambda=\frac 1 4 h^2=p^2.$$
So we have $\chi_m(D)=
\delta_m pi$, for some $\delta_m \in \{\pm 1\}$. Thus
\begin{align} \notag
\delta_m pi&=\chi_m(D)=\chi_m\left (\sum_{j,k,w} \varepsilon_{j,k,w} x^jy^kz^w    \right )\\
&= \sum_{j,k,w} \varepsilon_{j,k,w} \chi_m(x^jy^kz^w)=
\sum_{j,k,w} \varepsilon_{j,k,w} \zeta_p^{jm}i^w.\label{eq211}
\end{align}

We write (\ref {eq211})  as follows, splitting off the $j=0$ part:
\begin{align}
\label{eqn1}
\delta_m pi=\sum_{k,w} \varepsilon_{0,k,w}i^w+\sum_{j=1}^{p-1}\sum_{k,w} \varepsilon_{j,k,w}\zeta_p^{jm}i^w.
\end{align}

Summing over $0\le m<p$ we obtain
$$pi\sum_{m=0}^{p-1} \delta_m =p\sum_{k,w} \varepsilon_{0,k,w}i^w+\sum_{j=1}^{p-1}\sum_{k,w} \varepsilon_{j,k,w}i^w
\sum_{m=0}^{p-1} \zeta_p^{jm}=p\sum_{k,w} \varepsilon_{0,k,w}i^w.
$$

Thus
$$\sum_{k,w} \varepsilon_{0,k,w}i^w=i\sum_{m=0}^{p-1} \delta_m.
$$
Substituting this value of $\sum_{k,w} \varepsilon_{0,k,w}i^w$ in (\ref {eqn1}) we get
$
\delta_m pi=i\sum_{m=0}^{p-1} \delta_m +\sum_{j=1}^{p-1}\sum_{k,w} \varepsilon_{j,k,w}  \zeta_p^{jm}i^w,
$
or 
\begin{align} \label {eqd0} 0=\left (
-\delta_m p+\sum_{m=0}^{p-1} \delta_m \right ) +\sum_{j=1}^{p-1}\sum_{k,w} \varepsilon_{j,k,w}  \zeta_p^{jm}i^{w-1}.
\end{align}

Now (since $p$ is an odd prime) we have   $\mathbb Q(i) \cap \mathbb Q (\zeta_p)=\mathbb Q,$ and looking at the part of   (\ref {eqd0}) with $w=0,2,$   gives
\begin{align} \label {eqd10}
0= \sum_{j=1}^{p-1}\sum_{k} ( \varepsilon_{j,k,0}  -\varepsilon_{j,k,2} )\zeta_p^{jm}.
\end{align}

Now we use the fact that for $a_u \in \mathbb Q, 0\le u<p$, if $\sum_{u=0}^{p-1} a_u\zeta_p^u=0$, then $a_0=a_1=\cdots=a_{p-1}.$ 
Thus (\ref {eqd10})   gives  
$$0= \sum_{k} ( \varepsilon_{j,k,0}  -\varepsilon_{j,k,2} ),  \quad \text { for all }
j>0.$$ 

Now from  Theorem \ref {thm1}  we have  $|D \cap x^jH|=h/2=p$ for all $0<j<p$. Thus for $j \ne 0$  we have
$$
\sum_k \varepsilon_{j,k,0}+\varepsilon_{j,k,2}=p.
$$
But from the above we have:
$ \sum_{k} ( \varepsilon_{j,k,0}  -\varepsilon_{j,k,2} )=0,
$ and adding we get  
$$2\sum_k \varepsilon_{j,k,0}=p,
$$ a contradiction since $p$ is odd.
\qed\medskip

\noindent {\bf Case $G_{14}$}:   Let $f \in \mathbb N$  satisfy $f^2\equiv -1 \mod p$. Then
$$  G=G_{14}=\langle x,y,z|x^p,y^p,(x,y),z^4,x^z=x,y^z=y^f\rangle.$$
Thus  $G=\langle x\rangle \times \langle y,z\rangle$.
Now there are $p$ elements of order $2$, namely $(z^2)^{y^i}$ and they generate the normal subgroup 
 $\langle y,z^2\rangle$ of order $2p$, which must then be $H$.
  We also have the 
   linear character $\chi_m$ as in (\ref {eqn0}). Then, using this character, the  same argument as in the $G_{11}$ case  gives a contradiction. Thus $G_{14}$ is not a DRAD group.\qed

\section {The group $G_{15}$}

 Let $f \in \mathbb N$  satisfy $f^2\equiv -1 \mod p$. Since $p$ is prime, this implies that $p \equiv 1 \mod 4$. Then
$$  G=G_{15}=\langle x,y,z|x^p,y^p,(x,y),z^4,x^z=x^{-1},y^z=y^f\rangle.$$
Here $G=\langle x \rangle \rtimes \langle y,z\rangle$. Now there are $p$ involutions $(z^2)^{y^i},0\le i<p$, and they generate 
$\langle y,z^2\rangle \cong D_{2p}$, so we must have $H=\langle y,z^2\rangle \cong D_{2p}, H \triangleleft G$. 

This group does not succumb to any of the character-theoretic arguments that have worked in previous cases.

The proof will be by contradiction, so that we will assume that we have such a difference set $D\subset G$, where we  write 
$$
D=\sum_{g \in G} \varepsilon_g g \in \mathbb Z[G],
$$
 with $ \varepsilon_g \in \{0,1\}$ and $\varepsilon_g=0$ if $g \in H$. 
Then  we have 
\begin{align}
G=H+D+D^{-1}, \quad D \cap D^{-1}=D \cap H=\emptyset,\label{eqD}
\end{align}
 and 
$$D^{-1} =\sum_{g \in G} \varepsilon_g g^{-1}.
$$

The set-up for the proof of this case will then be as follows: we assume that 
$$D=\sum_{g \in G} \varepsilon_g g,
$$ where we can now think of the $\varepsilon_g$s as polynomial indeterminates generating the polynomial ring $R=\mathbb Z[\varepsilon_g]_{g \in G}$. Then $D$ and $D^{-1}$  will belong to the group  ring
$R[G]$. 

From (\ref{eqD})     we see that the  $\varepsilon_g$ satisfy
\begin{align}
\varepsilon_g^2=\varepsilon_g,\quad \varepsilon_g+\varepsilon_{g^{-1}}=1, \text { for } g \notin H, \quad \varepsilon_g=0 \text { for } g \in H.
\label{eqD1}
\end{align}

Let $\mathcal I$ denote the ideal of $\mathbb Z[\varepsilon_g]_{g \in G}$ generated 
by the relations  in (\ref{eqD1}) and $2\mathbb Z$.

We now consider the relations that we obtain from the difference set equation.
Let 
\begin{align}
E=DD^{-1}-(\lambda(G-1)+k) \in \mathbb Z[\varepsilon_g]_{g \in G}
\label{eqe} 
\end{align}
 and for $k \in G$ let $E_k$ denote the coefficient of $k$ in $E$. 
Then for $k \in G, k \ne 1,$ 
  we have
$E_k \in \mathbb Z[\varepsilon_g]_{g \in G}$. 

Since 
$$E=DD^{-1}  -(\lambda(G-1)+k)  =\sum_{g \in G}\varepsilon_g g\sum_{h \in G}\varepsilon_h h^{-1}-(\lambda(G-1)+k),
$$ we see that  for $k \in G, k \ne 1,$ we have (from (\ref{eqe})):
\begin{align}
E_k=  \sum_{h \in G} \varepsilon _{kh}\varepsilon_h-\lambda. \label{eqek}
\end{align}

  Using the relations in $G$  we see that a normal form for elements of $G$ is $g=x^iy^jz^k,$ where $ 0\le i,j<p,0\le k<4$.

We will need the following group-theoretic  facts about $G=\langle x,y,z\rangle$:
\begin{align*}
&(x^iy^j)^{-1}=x^{-i}y^{-j}, &(x^iy^jz^2)^{-1}=x^iy^{-j}z^2,\\&(x^iy^jz)^{-1}=x^{-i}y^{-fj}z^3,
&(x^iy^jz^3)^{-1}=x^{-i}y^{fj}z 
\end{align*}

From now on most of the equations that we write down will be considered to be in the quotient ring   $R=\mathbb Z[\varepsilon_g]_{g \in G}/\mathcal I$.  We will sometimes indicate this by writing $\mod \mathcal I$ at the end of the equation. 
  
  Thus, for example, we know from (\ref {eqek}) that, since $\lambda=p^2-p$ is even,  we have 
  $$E_k=
  \sum_{h \in G}\varepsilon_{kh}\varepsilon_h
   \mod \mathcal I
   $$
    for $k \in G, k \ne 1$.

Define subgroups
$$Y=\langle y\rangle, \quad N=\langle x,y,z^2\rangle, \quad A=\langle x,y\rangle,
$$
 so that $Y, N,A \triangleleft G$.   
 
   Let $p_2=(p-1)/2.$

Now for $k \in G$ summing the $E_g$ over the elements $g$ of the coset $Yk=kY$ we define:
\begin{align}
Z_k =\sum_{i=0}^{p-1} E_{y^ik}=
\sum_{i=0}^{p-1}  \sum_{h \in G} \varepsilon_{y^ikh}\varepsilon_h\in \mathbb Z[\varepsilon_g].\label{eqdef}
\end{align}

           For $g \in G$ we   define 
                  $$\mathcal Y(g)=\{y^ig:0\le i<p\}= \{gy^i:0\le i<p\}.
                  $$
                  We also define 
                  $$\Sigma_\mathcal Y(g)=
                  \sum_{u \in \mathcal Y(g)} \varepsilon_{u}
                  =\sum_{i=0}^{p-1} \varepsilon_{y^ig}.
                  $$ 
          
          Let $\mathcal E$ be the ideal of $\mathbb Z[\varepsilon_g]_{g \in G}$ generated by $\mathcal I$ and all the $ E_k, k\in G, k \ne 1$.
We note that if there is a DRAD difference set for $G$, then the quotient ring $R/\mathcal E$ would be non-trivial.
  Thus the following result will conclude the $G_{15}$ case:

\begin{thm} \label{thm111} As an equation in $\mathbb Z[\varepsilon_g]/\mathcal I$ we have
$$Z_x+Z_{xz^2}+Z_{x^{(p-1)/2}z^2}=1 \mod \mathcal I.$$

\end{thm}
   \noindent {\it Proof} 
Define $Z_1,\cdots,Z_6\in R$ as follows: 
\begin{align*}
&Z_1=   \sum_{i=0}^{p-1}    \sum_{h \in G \setminus N} \varepsilon_{y^ixh}\varepsilon_h;
\\
&Z_2=  \sum_{i=0}^{p-1}  \sum_{h \in G \setminus N} \varepsilon_{y^ixz^2h}\varepsilon_h;
\\
&Z_3=  \sum_{i=0}^{p-1}  \sum_{h \in G \setminus N} \varepsilon_{y^ix^{p_2}z^2h}\varepsilon_h;
\\
&Z_4=  \sum_{i=0}^{p-1}  \sum_{h \in N \setminus H} \varepsilon_{y^ixh}\varepsilon_h;
\\
&Z_5=  \sum_{i=0}^{p-1}  \sum_{h \in N \setminus H} \varepsilon_{y^ixz^2h}\varepsilon_h;\\
&Z_6=  \sum_{i=0}^{p-1}  \sum_{h \in N \setminus H} \varepsilon_{y^ix^{p_2}z^2h}\varepsilon_h.
\end{align*}

  Then we have: 
  $$Z_x=Z_1+Z_4, \quad Z_{xz^2}=Z_2+Z_5, \quad Z_{x^{p_2}z^2}=Z_3+Z_6,
  $$ 
  so that 
  $$Z_x+Z_{xz^2}+Z_{x^{p_2}z^2}=\sum_{i=1}^6 Z_i.
  $$

Thus we will now find $Z_1,\cdots,Z_6,$ and Theorem \ref {thm111} will follow by showing that $\sum_{i=1}^6 Z_i=1.$
Our first goal is simply stated:

\begin{lem}\label{lem1111} $Z_{1}=1,\,\,\, Z_{2}=0,\,\,\, Z_{3}=0.$
\end{lem}
   \noindent {\it Proof} 
We partition $G \setminus N$ into $P \cup P^{-1}$ where
\begin{align*}
&P=Az=\{x^ay^bz:0\le a,b<p\};  
\end{align*}
One can check that $\{P, P^{-1}\}$ is a partition of $G \setminus N$.

We consider the pairs $(y^ikh,h) \in (G\setminus N)^2$.   We wish to pair up a  pair
$(y^ikh,h)$    (see the definition of $Z_k$ in (\ref{eqdef}) above to motivate this) with another  pair of the form $(y^jkh_1,h_1)$, where
$(y^jkh_1,h_1)=((y^ikh)^{-1},h^{-1})$. This second  pair is  called the {\it dual} of $(y^ikh,h)$.
We now explain how to do this pairing for the situation where $h \in G \setminus N$ and $k \in \{x,xz^2,x^{p_2}z^2\}$, so that $k$ has the form $k=x^cz^{2d}$

Suppose that $h \in P \subset G \setminus N$, so that $h=x^ay^bz$. Then for $k=x^cz^{2d} $ we 
want the  pair dual to $(y^ikh,h)$ to 
be $((y^ikh)^{-1},h^{-1})$; but we want it to have the correct form i.e. 
to also 
look like $(y^j k  h_1,h_1)$  for some $j$ and some $h_1 \in G \setminus N$. 
Solving  
$((y^ikh)^{-1},h^{-1}) = (y^j k  h_1,h_1)$
we get 
 $h_1=h^{-1}$ and $(y^ikh)^{-1}= y^j k  h_1,$ and so
we need to find  $j$ so that  $h^{-1}k^{-1}y^{-i}\cdot h_1^{-1}k^{-1}   =y^j$ (if possible). Well from these equations we have
\begin{align}  \notag
h^{-1}k^{-1}y^{-i}\cdot h_1^{-1}k^{-1}    &=h^{-1}k^{-1}y^{-i}hk^{-1}   \\ \notag
&=z^{-1}x^{-a}y^{-b}\cdot x^{-c}z^{2d}\cdot y^{-i}\cdot x^ay^bz\cdot x^{-c}z^{2d}\\ \notag
&=z^{-1}y^{-b}\cdot x^{-c}z^{2d}\cdot y^{-i}\cdot y^bz\cdot x^{-c}z^{2d}\\ \notag
&=z^{-1}y^{-b}\cdot z^{2d}\cdot y^{-i}\cdot y^bz\cdot z^{2d}\\
&=y^{-fb}\cdot z^{2d}\cdot y^{-fi}\cdot y^{fb}\cdot z^{2d}\in Y\label{eqy}
\end{align}

Thus  for $h \in P$ and $k=x^cz^{2d},$ with the value of $j$ determined by (\ref{eqy}),   we pair $(y^ikh,h)$ with $((y^ikh)^{-1},h^{-1})=(y^jkh^{-1},h^{-1}).$

For $h \in P \subset G \setminus N, k\in \{x,xz^2,x^{p_2}z^2\}$, then, as described above,  we pair $(y^ikh,h) $ with $((y^ikh)^{-1},h^{-1})$.
The point here is that the part of the sum for $Z_k$ coming from the pairs  $(y^ikh,h) $ (for  $h \in P$ and summing over  $0\le i<p$) is 
\begin{align}
\label {eq21}
\sum_{i=0}^p \varepsilon_{y^ikh}  \varepsilon_h,
\end{align} 
while the sum coming from the dual pairs $((y^ikh)^{-1},h^{-1})$  (where $h^{-1} \in P^{-1}$) is
\begin{align}
\label {eq22}
\sum_{i=0}^{p-1} \varepsilon_{(y^ikh)^{-1}}  \varepsilon_{h^{-1}}= \sum_{i=0}^{p-1} (1-\varepsilon_{y^ikh})(1-  \varepsilon_{h}).
\end{align} 
Adding (\ref {eq21}) and (\ref {eq22})  we get (considering these elements modulo the ideal $\mathcal I$)
\begin{align}
\label {eq23}
\sum_{i=0}^{p-1} (1+\varepsilon_{y^ikh}+\varepsilon_h), \text { which is equal to }  1+\sum_{i=0}^{p-1} (\varepsilon_{y^ikh}+\varepsilon_h) .
\end{align}
Now we add (\ref {eq23}) over all $h=x^ay^bz  \in P$ to get     (since $p$ is odd) 
\begin{align}
\label {eq24}
\sum_{a,b=0}^{p-1} \left (1+\sum_{i=0}^{p-1} (\varepsilon_{y^ik x^ay^bz }+\varepsilon_{x^ay^bz }) \right) = 1+\sum_{a,b=0}^{p-1} \sum_{i=0}^{p-1} (\varepsilon_{y^ik x^ay^bz }+\varepsilon_{x^ay^bz })\mod \mathcal I.
\end{align}

Now  if $k =x$, then the sum (\ref {eq24}) is $Z_1$; but   since, with $k=x,$ we have 
$$\sum_{a,b=0}^{p-1}  (\varepsilon_{y^ik x^ay^bz }+\varepsilon_{x^ay^bz }) =\sum_{a,b=0}^{p-1}  (\varepsilon_{y^ix x^ay^bz }+\varepsilon_{x^ay^bz }) =0 \mod \mathcal I,
$$
 we see that $Z_{1}=1$.

On the other hand, if  $k=x^cz^2, c\in \{1,p_2\},$ then for (\ref {eq24}) we get 
 \begin{align*}
& 1+\sum_{i=0}^{p-1} \sum_{a,b=0}^{p-1} ( \varepsilon_{y^ix^cz^2x^ay^b z}+\varepsilon_{x^ay^bz})
\\&=
 1+\sum_{i=0}^{p-1} \sum_{a,b=0}^{p-1} ( \varepsilon_{x^ay^b z^3}+\varepsilon_{x^ay^bz})
= 1+\sum_{i=0}^{p-1} \sum_{a,b=0}^{p-1} 1=0 \mod \mathcal I.
\end{align*}
Thus $Z_{2}=Z_{3}=0$.
This concludes the proof of Lemma \ref {lem1111}.\qed\medskip

    Having found $Z_1, Z_2, Z_3$ we now show how to determine $Z_4, Z_5, Z_6$.

    For $g \in G$ recall that
    $
    \Sigma_\mathcal Y(g)=\sum_{i=0}^{p-1} \varepsilon_{y^ig}=\sum_{i=0}^{p-1} \varepsilon_{gy^i}.
   $
    
        \begin{lem}\label{lemUU}  (i)  For $g \in G\setminus H$ we have $\Sigma_\mathcal Y(g) + \Sigma_\mathcal Y(g^{-1})=1 \mod \mathcal I$.

 \noindent  (ii) For $h \in H$ we have $\Sigma_\mathcal Y(h)=0 \mod \mathcal I.$
 
  \noindent  (iii)  For $g \in G$ we have $\Sigma_\mathcal Y(g)\cdot \Sigma_\mathcal Y(g^{-1})=0 \mod \mathcal I.$

\end{lem}
\noindent {\it Proof}
(i) For $g \in G \setminus H$ we have, using the fact that $p$ is odd:
\begin{align*}
\Sigma_\mathcal Y(g) + \Sigma_\mathcal Y(g^{-1})&=
\sum_{i=0}^{p-1} \varepsilon_{y^ig}+\varepsilon_{y^ig^{-1}}
=\sum_{i=0}^{p-1} \varepsilon_{y^ig}+\varepsilon_{(y^ig)^{-1}}
=\sum_{i=0}^{p-1} 1 =1 \mod \mathcal I.
\end{align*} 
(ii) is clear since for $g \in H$ we have $y^ig \in H$, so that $\varepsilon_{y^ig}=0$.

  \noindent  
(iii) For $g \in H$ this follows from (ii), while for $g \in G \setminus H$    we have, using (i):
      \begin{align}
\Sigma_\mathcal Y(g)\cdot \Sigma_\mathcal Y(g^{-1})&=
\Sigma_\mathcal Y(g)\cdot (1+\Sigma_\mathcal Y(g)).
\label{eqwwp7}
\end{align}
Letting $E=\Sigma_\mathcal Y(g) $ we have $E^2=E \mod \mathcal I$ and so (\ref {eqwwp7}) is $E(1+E)=0 \mod \mathcal I$.  
\qed\medskip
              
We are now in a position  to prove:
         
         \begin{prop}\label{propZxzz}
         
        (i) $Z_4=  \Sigma_\mathcal Y(x)+ \Sigma_\mathcal Y(x^{p_2})+ \Sigma_\mathcal Y(xz^2)+ \Sigma_\mathcal Y(x^{p_2}z^2) \mod \mathcal I.$
         
  \noindent     (ii)    $Z_5=Z_{xz^2}=1+
\sum _{a=1}^{p-2}   ( \Sigma_\mathcal Y(x^{a}) +  \Sigma_\mathcal Y(x^{a+1}z^2) ) \mod \mathcal I;
    $
    
     \noindent     (iii) 
      \begin{align*}
      Z_6=Z_{x^{p_2}z^2}=&
      1+ \Sigma_\mathcal Y(x^{p_2}z^2) +
 \sum _{a=1}^{p_2}    \Sigma_\mathcal Y(x^{a}) +  \Sigma_\mathcal Y(x^{a+p_2}z^2))\\&
 +
\sum _{a=1}^{p_2}    \Sigma_\mathcal Y(x^{a+1+p_2}) +  \Sigma_\mathcal Y(x^{a}z^2) \mod \mathcal I.
\end{align*}
      \end{prop}
    \noindent {\it Proof} For the proof of Proposition \ref {propZxzz} (i), we have
  $Z_4=  \sum_{i=0}^{p-1}  \sum_{h \in N \setminus H} \varepsilon_{y^ixh}\varepsilon_h$ and $N\setminus H$ is a union of $Y$-cosets, where 
  $$
  (N\setminus H)/Y=\{x,x^2,\cdots,x^{p-1}\} \cup \{xz^2,x^2z^2,\cdots,x^{p-1}z^2\}
  $$
   is a set of coset representatives for $N \setminus H$. Thus we can write
\begin{align}
Z_4\notag &= \sum_{i=0}^{p-1}  \sum_{h \in N \setminus H} \varepsilon_{y^ixh}\varepsilon_h
= \sum_{i=0}^{p-1}\sum_{h \in (N\setminus H)/Y} \sum_{b=0}^{p-1} \varepsilon_{y^ihy^b}\varepsilon_{hy^b}
\\&= \sum_{i=0}^{p-1}  \sum_{b=0}^{p-1} \sum_{a=1}^{p-1} \varepsilon_{y^ixx^ay^b}\varepsilon_{x^ay^b} +\notag\varepsilon_{y^ixx^ay^bz^2}\varepsilon_{x^ay^bz^2} \\
&=\notag
\sum_{i=0}^{p-1}  \sum_{b=0}^{p-1} \sum_{a=1}^{p-1} \varepsilon_{y^ixx^a}\varepsilon_{x^ay^b}
+ \varepsilon_{y^ixx^az^2}\varepsilon_{x^ay^bz^2}\\
&= \sum_{a=1}^{p-1}\left ( \sum_{i=0}^{p-1}   \varepsilon_{y^ixx^a}    \sum_{b=0}^{p-1} \varepsilon_{x^ay^b}\right )  +\notag
\sum_{a=1}^{p-1}\left ( \sum_{i=0}^{p-1}   \varepsilon_{y^ixx^az^2}    \sum_{b=0}^{p-1} \varepsilon_{x^ay^bz^2}\right )  \\
&=
 \sum_{a=1}^{p-1}  \Sigma_\mathcal Y(x^{a+1}) \Sigma_\mathcal Y(x^a)  +   \sum_{a=1}^{p-1}  \Sigma_\mathcal Y(x^{a+1}z^2) \Sigma_\mathcal Y(x^az^2).  \label{eqwrt}
\end{align}
 Let the   two sums in (\ref {eqwrt})  be denoted by $Z_{41},Z_{42}$. Then using Lemma \ref {lemUU} and the fact that $p_2$ is even, we have:
 \begin{align}
Z_{41}&= \sum_{a=1}^{p-1}  \Sigma_\mathcal Y(x^{a+1}) \Sigma_\mathcal Y(x^a) \notag \\
&=\notag  \sum_{a=1}^{p_2}  \Sigma_\mathcal Y(x^{a+1}) \Sigma_\mathcal Y(x^a) +
 \sum_{a=p_2+1}^{p-1}  \Sigma_\mathcal Y(x^{a+1}) \Sigma_\mathcal Y(x^a) \notag\\
 &=
\left (  \sum_{a=1}^{p_2-1}  \Sigma_\mathcal Y(x^{a+1}) \Sigma_\mathcal Y(x^a) \right ) +  \Sigma_\mathcal Y(x^{p_2+1}) \Sigma_\mathcal Y(x^{p_2}) 
+ \sum_{a=p_2+1}^{p-1}  \Sigma_\mathcal Y(x^{a+1}) \Sigma_\mathcal Y(x^a) \notag \\
 &=
\left (  \sum_{a=1}^{p_2-1}  \Sigma_\mathcal Y(x^{a+1}) \Sigma_\mathcal Y(x^a) \right ) +  \Sigma_\mathcal Y(x^{p_2+1}) \Sigma_\mathcal Y(x^{p_2}) 
+ \sum_{a=1}^{p_2-1}  \Sigma_\mathcal Y(x^{a+1+p_2}) \Sigma_\mathcal Y(x^{a+p_2}) \notag \\
&=
\left (  \sum_{a=1}^{p_2-1}  \Sigma_\mathcal Y(x^{a+1}) \Sigma_\mathcal Y(x^a) \right ) +  \Sigma_\mathcal Y(x^{p_2+1}) \Sigma_\mathcal Y(x^{p_2}) 
+ \sum_{a=1}^{p_2-1}  \Sigma_\mathcal Y(x^{-a}) \Sigma_\mathcal Y(x^{-a-1}) \notag \\
 &=
  \sum_{a=1}^{p_2-1}  \Sigma_\mathcal Y(x^{a+1}) \Sigma_\mathcal Y(x^a) + 0 
+ \sum_{a=1}^{p_2-1} (1+ \Sigma_\mathcal Y(x^{a})) (1+\Sigma_\mathcal Y(x^{a+1})) \notag \\
&=
\sum_{a=1}^{p_2-1} (1+ \Sigma_\mathcal Y(x^{a})+\Sigma_\mathcal Y(x^{a+1})) \notag \\
&=1+
\sum_{a=1}^{p_2-1}  \Sigma_\mathcal Y(x^{a})+\Sigma_\mathcal Y(x^{a+1}).\notag
\end{align}
For $Z_{42} $ we similarly obtain
 \begin{align*} 
Z_{42}=1+
\sum_{a=1}^{p_2-1}  \Sigma_\mathcal Y(x^{a}z^2)+\Sigma_\mathcal Y(x^{a+1}z^2).
\end{align*}
Thus 
\begin{align*}
Z_4&=Z_{41}+Z_{42}\\
&=\sum_{a=1}^{p_2-1} ( \Sigma_\mathcal Y(x^{a})+\Sigma_\mathcal Y(x^{a+1}))
+
\sum_{a=1}^{p_2-1} ( \Sigma_\mathcal Y(x^{a}z^2)+\Sigma_\mathcal Y(x^{a+1}z^2))\\
&= \Sigma_\mathcal Y(x)+ \Sigma_\mathcal Y(x^{p_2})+ \Sigma_\mathcal Y(xz^2)+ \Sigma_\mathcal Y(x^{p_2}z^2),
\end{align*}
as required for Proposition \ref {propZxzz} (i).

\medskip

For $Z_5$ we again sum over cosets of $Y$, so that we have
 \begin{align}    \notag
Z_5&=\sum_{i=0}^{p-1}\sum _{h \in N \setminus H} \varepsilon_{y^ixz^2h} \varepsilon_h\\
&=\sum_{i,b=0}^{p-1}\sum _{a=1}^{p-1}  \varepsilon_{y^ixz^2x^ay^b} \varepsilon_{x^ay^b}
+  \notag
\sum_{i,b=0}^{p-1}\sum _{a=1}^{p-1}  \varepsilon_{y^ixz^2x^ay^bz^2} \varepsilon_{x^ay^bz^2}\\
&=\sum_{i,b=0}^{p-1}\sum _{a=1}^{p-1}  \varepsilon_{y^ixz^2x^a} \varepsilon_{x^ay^b}  \notag
+
\sum_{i,b=0}^{p-1}\sum _{a=1}^{p-1}  \varepsilon_{y^ixz^2x^az^2} \varepsilon_{x^ay^bz^2}\\
&=
\sum _{a=1}^{p-1}  \Sigma_\mathcal Y(x^{a+1}z^2)  \Sigma_\mathcal Y(x^a)
+
\sum _{a=1}^{p-1}    \Sigma_\mathcal Y(x^{a+1})  \Sigma_\mathcal Y(x^az^2).\label{eq1222}
\end{align}
Now when  $a=p-1$ the corresponding terms of each sum of  (\ref {eq1222}) are zero, while all other terms are non-zero.
 Thus we have
\begin{align*}
Z_5&=
\sum _{a=1}^{p-2}  \Sigma_\mathcal Y(x^{a+1}z^2)  \Sigma_\mathcal Y(x^a)
+  
\sum _{a=1}^{p-2}   (1+ \Sigma_\mathcal Y(x^{-a-1})) (1+ \Sigma_\mathcal Y(x^{-a}z^2))\\
&=
\sum _{a=1}^{p-2}  \Sigma_\mathcal Y(x^{a+1}z^2)  \Sigma_\mathcal Y(x^a)\\
&\qquad 
+  
\sum _{a=1}^{p-2}   (1+ \Sigma_\mathcal Y(x^{-a-1}) +  \Sigma_\mathcal Y(x^{-a}z^2) +  
 \Sigma_\mathcal Y(x^{p-a}z^2)  \Sigma_\mathcal Y(x^{p-a-1})  )\\
 &=
\sum _{a=1}^{p-2}  \Sigma_\mathcal Y(x^{a+1}z^2)  \Sigma_\mathcal Y(x^a)\\&\qquad 
+  1+
\sum _{a=1}^{p-2}   ( \Sigma_\mathcal Y(x^{-a-1}) +  \Sigma_\mathcal Y(x^{-a}z^2) +  
 \Sigma_\mathcal Y(x^{p-a}z^2)  \Sigma_\mathcal Y(x^{p-a-1})  )
\end{align*}

Now one checks that 
$\sum _{a=1}^{p-2}  \Sigma_\mathcal Y(x^{a+1}z^2)  \Sigma_\mathcal Y(x^a) =
\sum _{a=1}^{p-2}   
 \Sigma_\mathcal Y(x^{p-a}z^2)  \Sigma_\mathcal Y(x^{p-a-1}) ,$ 
so that we have
\begin{align*}
Z_5&=1+\sum _{a=1}^{p-2}   ( \Sigma_\mathcal Y(x^{-a-1}) +  \Sigma_\mathcal Y(x^{-a}z^2) )
=  1+
\sum _{a=1}^{p-2}   ( \Sigma_\mathcal Y(x^{a}) +  \Sigma_\mathcal Y(x^{a+1}z^2) ).
\end{align*}
This gives Proposition \ref {propZxzz} (ii).

\medskip

For $Z_6$ we again sum over cosets of $Y$, so that we have
 \begin{align}   \notag 
Z_6&=\sum_{i=0}^{p-1}\sum _{h \in N \setminus H} \varepsilon_{y^ix^{p_2}z^2h} \varepsilon_h\\
&=\sum_{i,b=0}^{p-1}\sum _{a=1}^{p-1}  \varepsilon_{y^ix^{p_2}z^2x^ay^b} \varepsilon_{x^ay^b}
+   \notag 
\sum_{i,b=0}^{p-1}\sum _{a=1}^{p-1}  \varepsilon_{y^ix^{p_2}z^2x^ay^bz^2} \varepsilon_{x^ay^bz^2}\\
&=\sum_{i,b=0}^{p-1}\sum _{a=1}^{p-1}  \varepsilon_{y^ix^{p_2}z^2x^a} \varepsilon_{x^ay^b}
+   \notag 
\sum_{i,b=0}^{p-1}\sum _{a=1}^{p-1}  \varepsilon_{y^ix^{p_2}z^2x^az^2} \varepsilon_{x^ay^bz^2}\\
&=  \label{eqp1rst}
\sum _{a=1}^{p-1}  \Sigma_\mathcal Y(x^{a+p_2}z^2)  \Sigma_\mathcal Y(x^a)
+
\sum _{a=1}^{p-1}    \Sigma_\mathcal Y(x^{a+p_2})  \Sigma_\mathcal Y(x^az^2)
\end{align}
Now in the sums of  (\ref {eqp1rst}) the only    term of  (\ref {eqp1rst}) having the form  $ \Sigma_\mathcal Y(*) \Sigma_\mathcal Y(*)$   that is zero,
occurs  when $a=p_2+1$, there being two such occurrences in  (\ref {eqp1rst}).
So  (\ref {eqp1rst}) is equal to 
 \begin{align}   \notag 
&   \notag 
\sum _{a=1,a \ne p_2+1}^{p-1}  \Sigma_\mathcal Y(x^{a+p_2}z^2)  \Sigma_\mathcal Y(x^a)
+  
  (1+ \Sigma_\mathcal Y(x^{-a-p_2})) (1+ \Sigma_\mathcal Y(x^{-a}z^2))\\
&=   \notag 
\sum _{a=1,a \ne p_2+1}^{p-1}  \Sigma_\mathcal Y(x^{a+p_2}z^2)  \Sigma_\mathcal Y(x^a)
+  \\&\qquad 
   (1+ \Sigma_\mathcal Y(x^{-a-p_2}) +  \Sigma_\mathcal Y(x^{-a}z^2) +     \notag 
 \Sigma_\mathcal Y(x^{p-a}z^2)  \Sigma_\mathcal Y(x^{p-a-p_2})  )\\
 &=1+
\sum _{a=1,a\ne p_2+1}^{p-1}  \Sigma_\mathcal Y(x^{a+p_2}z^2)  \Sigma_\mathcal Y(x^a)    \notag 
+  \\
&\qquad 
\sum _{a=1,a \ne p_2+1}^{p-1}   ( \Sigma_\mathcal Y(x^{-a-p_2}) +  \Sigma_\mathcal Y(x^{-a}z^2) +  
 \Sigma_\mathcal Y(x^{p-a}z^2)  \Sigma_\mathcal Y(x^{p-a-p_2})  ).\label{eqd12}
\end{align}
We split (\ref{eqd12})  into three parts,   so that $Z_6=1+Z_{61}+Z_{62}$. Taking $a=1,\cdots,p_2$ in (\ref {eqd12}) will determine
$Z_{61}$:
\begin{align} Z_{61}=   & \notag 
\sum _{a=1}^{p_2}  \Sigma_\mathcal Y(x^{a+p_2}z^2)  \Sigma_\mathcal Y(x^a)
+  \\
&\qquad 
\sum _{a=1}^{p_2}   ( \Sigma_\mathcal Y(x^{-a-p_2}) +  \Sigma_\mathcal Y(x^{-a}z^2) +  
 \Sigma_\mathcal Y(x^{p-a}z^2)  \Sigma_\mathcal Y(x^{p-a-p_2})  ),\label{eqd2}
\end{align}
where we note that, since for any function $f$ we have $\sum_{a=1} ^{p_2} f(a)=\sum_{a=1}^{p_2} f(1+p_2-a)$, the degree $2$ part of (\ref {eqd2})   is equal to
\begin{align*}
\sum _{a=1}^{p_2} & \Sigma_\mathcal Y(x^{a+p_2}z^2)  \Sigma_\mathcal Y(x^a)
+  
\sum _{a=1}^{p_2} 
 \Sigma_\mathcal Y(x^{p-a}z^2)  \Sigma_\mathcal Y(x^{p-a-p_2})  \\
 &=
 \sum _{a=1}^{p_2}  \Sigma_\mathcal Y(x^{a+p_2}z^2)  \Sigma_\mathcal Y(x^a)
+  
\sum _{a=1}^{p_2} 
 \Sigma_\mathcal Y(x^{(p-(p_2+1-a))}z^2)  \Sigma_\mathcal Y(x^{p_2+1-(p_2+1-a)})  \\
 &=
 \sum _{a=1}^{p_2}  \Sigma_\mathcal Y(x^{a+p_2}z^2)  \Sigma_\mathcal Y(x^a)
+  
\sum _{a=1}^{p_2} 
 \Sigma_\mathcal Y(x^{a+p_2)}z^2)  \Sigma_\mathcal Y(x^{a)})   =0.
\end{align*}
With this we now have
\begin{align}
Z_{61}&=\sum _{a=1}^{p_2}    \Sigma_\mathcal Y(x^{-a-p_2}) +  \Sigma_\mathcal Y(x^{-a}z^2)\notag\\
&=\notag 
\sum _{a=1}^{p_2}    \Sigma_\mathcal Y(x^{-(1+p_2-a)-p_2}) +  \Sigma_\mathcal Y(x^{-(1+p_2-a)}z^2)\\
&=
\sum _{a=1}^{p_2}    \Sigma_\mathcal Y(x^{a}) +  \Sigma_\mathcal Y(x^{a+p_2}z^2).
\label{eq613}
\end{align}

Taking $a=p_2+2,\cdots ,p-1$ in (\ref {eqd12}) gives $Z_{62}$:
\begin{align}
Z_{62}&=
    \notag 
\sum _{a=2+p_2}^{p-1}  \Sigma_\mathcal Y(x^{a+p_2}z^2)  \Sigma_\mathcal Y(x^a)
+  \\
&\qquad 
\sum _{a=2+p_2}^{p-1}   ( \Sigma_\mathcal Y(x^{-a-p_2}) +  \Sigma_\mathcal Y(x^{-a}z^2) +  
 \Sigma_\mathcal Y(x^{p-a}z^2)  \Sigma_\mathcal Y(x^{p-a-p_2})  ) \label{eqd2}.
\end{align}
Again we look at the  degree $2$ part of (\ref{eqd2}):
\begin{align}&
\sum _{a=2+p_2}^{p-1}  \Sigma_\mathcal Y(x^{-a-p_2}z^2)  \Sigma_\mathcal Y(x^a)
+  
 \Sigma_\mathcal Y(x^{p-a}z^2)  \Sigma_\mathcal Y(x^{p-a-p_2}),\label{eqp122}  
  \end{align}
  but it is easy to see that the terms in each of the sums of (\ref {eqp122}) 
 are the same only listed in reverse order in the second sum, so that 
  (\ref {eqp122}) is
 zero. Thus the degree $2$ part of  (\ref{eqd2}) is zero, so that
 \begin{align}
 Z_{62}&=\notag 
\sum _{a=2+p_2}^{p-1}    \Sigma_\mathcal Y(x^{-a-p_2}) +  \Sigma_\mathcal Y(x^{-a}z^2)\\
&=\notag 
\sum _{a=2}^{p_2}    \Sigma_\mathcal Y(x^{-a-2p_2}) +  \Sigma_\mathcal Y(x^{-a-p_2}z^2)\\
&=
\sum _{a=2}^{p_2}    \Sigma_\mathcal Y(x^{1-a}) +  \Sigma_\mathcal Y(x^{1+p_2-a}z^2)  \notag  \\
&= \Sigma_\mathcal Y(x^{p_2}z^2) +
\sum _{a=1}^{p_2}    \Sigma_\mathcal Y(x^{1-a}) +  \Sigma_\mathcal Y(x^{1+p_2-a}z^2)  \notag  \\
&= \Sigma_\mathcal Y(x^{p_2}z^2) +
\sum _{a=1}^{p_2}    \Sigma_\mathcal Y(x^{a+1+p_2}) +  \Sigma_\mathcal Y(x^{a}z^2). \label{eq623}
  \end{align} 
 From (\ref {eq613}) and (\ref {eq623}) we get 
 \begin{align*}
 Z_6=1+Z_{61}+Z_{62}=&
 1+ \Sigma_\mathcal Y(x^{p_2}z^2) +
 \sum _{a=1}^{p_2}    \Sigma_\mathcal Y(x^{a}) +  \Sigma_\mathcal Y(x^{a+p_2}z^2)\\&
 +
\sum _{a=1}^{p_2}    \Sigma_\mathcal Y(x^{a+1+p_2}) +  \Sigma_\mathcal Y(x^{a}z^2),
\end{align*}
which completes the proof of Proposition  \ref {propZxzz}.\qed\medskip

\begin{cor} \label{cor1} $Z_4+Z_5+Z_6=0.$
\end{cor} 
\noindent{\it Proof}
 From Proposition  \ref {propZxzz} we have
 \begin{align}&Z_4+Z_5+Z_6=  \notag 
  \Sigma_\mathcal Y(x)+ \Sigma_\mathcal Y(x^{p_2})+ \Sigma_\mathcal Y(xz^2)\\&
  \qquad   \notag 
  +
        1+
\sum _{a=1}^{p-2}   ( \Sigma_\mathcal Y(x^{a}) +  \Sigma_\mathcal Y(x^{a+1}z^2) )
\\
 &  \notag 
 +1   +  \Sigma_\mathcal Y(x^{p_2}z^2)   +   \sum _{a=1}^{p_2}  (  \Sigma_\mathcal Y(x^{a}) +  \Sigma_\mathcal Y(x^{a+p_2}z^2))
 +
\sum _{a=1}^{p_2}  (  \Sigma_\mathcal Y(x^{a+1+p_2}) +  \Sigma_\mathcal Y(x^{a}z^2))\\
&=
  \Sigma_\mathcal Y(x)+ \Sigma_\mathcal Y(x^{p_2})+ \Sigma_\mathcal Y(xz^2)+ \Sigma_\mathcal Y(x^{p_2}z^2)
  +  \notag 
\sum _{a=1}^{p-2}   ( \Sigma_\mathcal Y(x^{a}) +  \Sigma_\mathcal Y(x^{a+1}z^2) )
\\
 &\qquad 
    +   \sum _{a=1}^{p_2}  (  \Sigma_\mathcal Y(x^{a}) +  \Sigma_\mathcal Y(x^{a+p_2}z^2) )
 +
\sum _{a=1}^{p_2}  (  \Sigma_\mathcal Y(x^{a+1+p_2}) +  \Sigma_\mathcal Y(x^{a}z^2)).\label{eqz456}
\end{align}
 Now one checks that
 $$
  \Sigma_\mathcal Y(x)+ \Sigma_\mathcal Y(x^{p_2})+
  \sum _{a=1}^{p-2}    \Sigma_\mathcal Y(x^{a})   
   +   \sum _{a=1}^{p_2}   \Sigma_\mathcal Y(x^{a})  +
 \sum _{a=1}^{p_2}    \Sigma_\mathcal Y(x^{a+1+p_2}) =0, $$
 so that (\ref{eqz456})  now gives
 \begin{align*}&Z_4+Z_5+Z_6=
 \Sigma_\mathcal Y(xz^2)+
 \sum _{a=1}^{p-2}    \Sigma_\mathcal Y(x^{a+1}z^2)   +
 \sum _{a=1}^{p_2}    \Sigma_\mathcal Y(x^{a}z^2)   +
  \sum _{a=1}^{p_2}    \Sigma_\mathcal Y(x^{a+1+p_2}z^2),
 \end{align*}
 which one can similarly see is equal to zero.
 
 Thus we have $Z_4+Z_5+Z_6=0$, concluding the proof of Corollary \ref {cor1}.\qed\medskip
 
So  from Lemma \ref {lem1111} and Corollary \ref {cor1} we have 
$$Z_x+Z_{xz^2}+Z_{x^{p_2}z^2}=\sum_{i=1}^6 Z_i=(Z_1+Z_2+Z_3)+(Z_4+Z_5+Z_6)=1+0=1,
$$
 concluding the proof of Theorem
 \ref {thm111}.\qed\medskip 
 
 Theorem \ref{thm2} follows as we have now considered (elliminated) each of the six groups not covered by the paper of Iiams \cite {ii}.\qed


\begin{thebibliography}{HJ}
 
 \bibitem{as}  
 AbuGhneim, Omar A.; Smith, Ken W.
\emph{Tightening Turyn's bound for Hadamard difference sets.}
J. Algebraic Combin. 27 (2008), no. 2, 187--203. 

\bibitem{Ma}  W. Bosma and J. Cannon,  MAGMA  (University of
Sydney, Sydney, 1994).



\bibitem{co}
H. Cohen, \emph{A Course in Computational Algebraic Number Theory}, GTM, vol. 138, Springer, 1996. 



\bibitem{DJ}
J. Davis and J. Jedwab, \emph{A survey of Hadamard difference sets}, HPL-94-14,
HP Laboratories, Bristol 1994.


\bibitem{dp}
Davis, James A.; Polhill, John \emph{Difference set constructions of DRADs and association schemes.} J. Combin. Theory Ser. A 117 (2010), no. 5, 598--605.











\bibitem{DS1} Courtney Hoagland, Stephen 
P. Humphries, Nathan Nicholson, Seth Poulsen
\emph{Difference Sets Disjoint from a Subgroup},
Graphs and Combinatorics (2019) 
35, 579--597 https://doi.org/10.1007/s00373-019-02017-2ORIGINAL 






\bibitem{ii}
Iiams, J., \emph{On difference sets in groups of order $4p^2$}, Journal of Comb. Theory A (1996) pp 256--276.
 
 
 \bibitem{ito1}  Ito, Noboru; Raposa, Blessilda P. \emph{Nearly triply regular DRADs of RH type.} Graphs Combin. 8 (1992), no. 2, 143--153. 
 

 
 \bibitem{ito2}
 Ito, Noboru \emph{Automorphism groups of DRADs.} Group theory (Singapore, 1987), 151--170, de Gruyter, Berlin, (1989).



 \bibitem{jed} 
 J. Jedwab, \emph {Perfect Arrays, Barker Arrays, and Difference Sets}, Ph.D. thesis, University of London, London, England (1991).

\bibitem{km} 
Kesava Menon, P. \emph{On difference sets whose parameters satisfy a certain relation.}  Proc. Amer. Math. Soc. 13, (1962) 739--745. 

 
  \bibitem{kib} 
 Kibler, Robert E.
\emph{A summary of noncyclic difference sets, $k<20$}. 
J. Combinatorial Theory Ser. A 25 (1978), no. 1, 62--67. 
 
 

 \bibitem{kra} 
 R. Kraemer, \emph{A result on Hadamard difference sets}, J. Combin. Theory (A), Vol. 63 (1993) pp. 1--10.

  \bibitem{mcf} 
 McFarland, Robert L.,
\emph{Difference sets in abelian groups of order $4p^2$}
Mitt. Math. Sem. Giessen No. 192 (1989), i--iv, 1--70.

  
 
 

 

 
 \bibitem{sm} 
 Smith, Ken W.
\emph{Non-abelian Hadamard difference sets.}
J. Combin. Theory Ser. A 70 (1995), no. 1, 144--156. 


 \bibitem{tur} 
 R. J. Turyn, \emph {Character sums and difference sets}. Pacific J. Math., Vol. 15 (1965) pp. 319--346.
 
   \bibitem{wan} 
 Wan, Z., \emph {Difference sets in groups of order $4p^4$}.
  Beijing Dexue Xuebao Ziran Kexue Ban 36(3), 331--341 (2000).
 
 
  \bibitem{web} 
 Webster, Jordan D. \emph{Reversible difference sets with rational idempotents.} Arab. J. Math. (Springer) 2 (2013), no. 1, 103--114.
 
 




 \end{thebibliography}
  \end{document}